\documentstyle{article}
\begin{document}

\def\dinf{\displaystyle\inf}
\def\dsup{\displaystyle\sup}
\def \theequation {\arabic {equation}}
\setcounter {equation}{0}

\vspace*{1cm}
\begin{center}
{{{\LARGE {\bf Robust Strictly Positive Real Synthesis for Convex
Combination of Sixth-Order Polynomials}\footnote{Supported by
National Natural Science Foundation of China and Natural Science
Foundation of National Laboratory of Intelligent Technology and
System of Tsinghua University.}}}}
\end{center}


\vskip 0.6cm \centerline{ Wensheng YU} \vskip 6pt

\centerline{\small {\it Institute of Automation, Chinese Academy
of Sciences}}

\centerline{\small\it {Beijing 100080, CHINA. E-mail:
wensheng.yu@mail.ia.ac.cn}}

\vskip 0.5cm \centerline{ Long WANG} \vskip 6pt

\centerline{\small {\it Center for Systems and Control,}}

\centerline{\small {\it Department of Mechanics and Engineering
Science, Peking University,}}

\centerline{\small\it {Beijing  100871, CHINA. E-mail:
longwang@mech.pku.edu.cn}}

\vskip 0.6cm \noindent { {{Abstract\quad}}} {{{For the two
sixth-order polynomials $a(s)$ and $b(s),$ Hurwitz stability of
their convex combination is necessary and sufficient for the
existence of a polynomial $c(s)$ such that $c(s)/a(s)$ and
$c(s)/b(s)$ are both strictly positive real. Our reasoning method
is constructive, and is insightful and helpful in solving the
general robust strictly positive real synthesis problem.}}}

\noindent
{{ {Keywords\quad}}}
{{{Robust Stability, Strict Positive Realness, Robust Analysis and Synthesis}}

\strut

{ { The strict positive realness (SPR) of transfer functions is an
important performance specification, and plays a critical role in
various fields such as absolute stability/hyperstability theory
\cite{Kal63,Pop73}, passivity analysis \cite{DV75}, quadratic
optimal control \cite{AM70} and adaptive system theory
\cite{Lan79}. In recent years, stimulated by the parametrization
method in robust stability analysis \cite{Bar94,BCK95}, the study
of robust strictly positive real systems has received much
attention, and great progress has beem made
\cite{ADK90,BZ93,BTV01,CDB91,
DB87,HHX89,HHX90,MA98,MP01,PD97,WH91,
WY99,WY00,WY01,WY01a,Yu98,YH99,YW99,YW00,YW01,YW01a,YW01b,YWT99}.
However, most results belong to the category of robust SPR
analysis. Valuable results in robust SPR synthesis are rare. The
following fundamental problem is still open
\cite{ADK90,BZ93,DB87,HHX89,HHX90,MA98,MP01,PD97,WH91,
WY01a,YH99,YW00,YW01,YW01a,YW01b}:

{\it Suppose $a(s)$ and $b(s)$ are two $n$-$th$ order Hurwitz
polynomials, does there exist, and how to find a (fixed)
polynomial $c(s)$ such that $c(s)/a(s)$ and $c(s)/b(s)$ are both
SPR?}

By the definition of SPR, it is easy to know that the Hurwitz
stability of the convex combination of $a(s)$ and $b(s)$ is
necessary for the existence of polynomial $c(s)$ such that
$c(s)/a(s)$ and $c(s)/b(s)$ are both SPR. In
\cite{HHX89,HHX90,PD97}, it was proved that, if $a(s)$ and $b(s)$
have the same even (or odd) parts, such a polynomial $c(s)$ always
exists; In \cite{ADK90, HHX89,HHX90,MA98,
WY99,WY00,WY01a,Yu98,YW00,YW01}, it was proved that, if $n\leq 4$
and $a(s),b(s)\in K$ ($K$ is a stable interval polynomial set),
such a polynomial $c(s)$ always exists; Recent results show that
\cite{WY99,WY00,WY01a,Yu98,YH99,YW00,YW01a,YW01b}, if $n\leq 5$
and $a(s)$ and $b(s)$ are the two endpoints of the convex
combination of stable polynomials, such a polynomial $c(s)$ always
exists. Some sufficient condition for robust SPR synthesis are
presented in \cite{ADK90,BZ93,DB87,MA98,WY99,WY00,Yu98},
especially, the design method in \cite{WY99,WY00} is numerically
efficient for high-order polynomial segments and interval
polynomials, and the derived conditions are necessary and
sufficient for low-order polynomial segments and interval
polynomials.

This paper shows that, for the two sixth-order polynomials $a(s)$
and $b(s),$ Hurwitz stability of their convex combination is
necessary and sufficient for the existence of a polynomial $c(s)$
such that $c(s)/a(s)$ and $c(s)/b(s)$ are both SPR. This also
shows that the conditions given in \cite{WY99,WY00} are also
necessary and sufficient, and the open problem above has a
positive answer for the case of sixth-order polynomial segment.
Our reasoning method is constructive, and is useful in solving the
general robust SPR synthesis problem.

In this paper, $ P^n$ stands for the set of $n$-$th$ order polynomials of $s$ with real coefficients,
$ R$ stands for the field of real numbers, $ \partial (p)$ stands for the order of polynomial
$ p(\cdot )$, and $H^n\subset P^n$ stands for the set of $n$-$th$ order Hurwitz stable polynomials
with real coefficients.

In the following definition, $ p(\cdot )\in P^m, q(\cdot )\in P^n, f(s)=p(s)/q(s)$ is
a rational function.

{\bf  Definition 1} \cite{YW99} \ \ $ f(s)$ is said to be strictly
positive real(SPR), if

(i) $ \partial (p)=\partial (q);$

(ii) $f(s)$ is analytic in $ \mbox{Re}  [s]\geq 0,$  (namely, $ q(\cdot )\in H^n$ );

(iii) $ \mbox{Re}  [f(j\omega )]>0,\ \ \ \forall \omega \in R.$

If $f(s)=p(s)/q(s)$ is proper, it is easy to get the following property:

{\bf Property 1} \cite{CDB91} \ \ If $f(s)=p(s)/q(s)$ is a proper
rational function, $q(s)\in H^n,$ and $\forall \omega \in R,
\mbox{Re}  [f(j\omega )]>0,$ then $p(s)\in H^n\cup H^{n-1}.$

The following theorem is the main result of this paper:

{\bf  Theorem 1} \ \ Suppose $
a(s)=s^6+a_1s^5+a_2s^4+a_3s^3+a_4s^2+a_5s+a_6\in H^6,
b(s)=s^6+b_1s^5+b_2s^4 +b_3s^3+b_4s^2+b_5s+b_6\in H^6,$ the
necessary and sufficient condition for the existence of a
polynomial $ c(s)$ such that $ c(s)/a(s)$ and $ c(s)/b(s)$ are
both Strictly Positive Real is
$$ \lambda b(s)+(1-\lambda )a(s)\in H^6,\lambda \in [0,1].$$

Since SPR transfer functions enjoy convexity property, by Property
1, we can easily get the necessary part of the theorem.

To prove sufficiency, we must first introduce some lemmas.

{\bf Lemma {\bf 1}}\ \ Suppose $
s^6+a_1s^5+a_2s^4+a_3s^3+a_4s^2+a_5s+a_6\in H^6,$ then the
following quadratic curve is an  ellipse in  the first quadrant of
the $x$-$y$-$z$-$p$ space:
\begin{center}
$\left\{
\begin{array}{l}
(a_2x+z-a_1y-a_3)^2-4(a_1-x)(a_5+a_3y+a_1p-a_2z-a_4x)=0 \\
a_6x+a_4z-a_3p-a_5y=0 \\
a_5p-a_6z=0
\end{array}
\right. $
\end{center}
and this ellipse is tangent with the line
\begin{center}
$\left\{
\begin{array}{l}
a_6x+a_4z-a_3p-a_5y=0 \\
a_5p-a_6z=0 \\
a_1-x=0
\end{array}
\right. $
\end{center}
at
\begin{center}
$\left\{
\begin{array}{l}
x=a_1, \\
y=\frac{a_5a_6a_1-a_5a_4a_2a_1+a_5a_4a_3+a_3a_6a_2a_1-a_3^2a_6}{%
-a_5a_4a_1+a_5^2+a_3a_6a_1}, \\
z=-a_5\frac{a_5a_2a_1-a_6a_1^2-a_5a_3}{-a_5a_4a_1+a_5^2+a_3a_6a_1}, \\
p=-a_6\frac{a_5a_2a_1-a_6a_1^2-a_5a_3}{-a_5a_4a_1+a_5^2+a_3a_6a_1},
\end{array}
\right. $
\end{center}
tangent with the line
\begin{center}
$\left\{
\begin{array}{l}
a_6x+a_4z-a_3p-a_5y=0 \\
a_5p-a_6z=0 \\
a_5+a_3y+a_1p-a_2z-a_4x=0
\end{array}
\right. $
\end{center}
at
\begin{center}
$\left\{
\begin{array}{l}
x=\frac{-a_5^3+a_5^2a_4a_1+a_3^3a_6+a_5^2a_2a_3-a_5a_3^2a_4-2a_5a_3a_6a_1}{%
a_5^2a_2^2-a_5^2a_4+a_5a_4^2a_1-2a_5a_6a_1a_2-a_3a_6a_1a_4+a_6^2a_1^2+a_5a_3a_6-a_5a_3a_4a_2+a_3^2a_6a_2%
}, \\
y=\frac{a_6a_4a_3^2-a_3a_6^2a_1-a_3a_5a_4^2-a_5^2a_6+a_5^2a_4a_2}{%
a_5^2a_2^2-a_5^2a_4+a_5a_4^2a_1-2a_5a_6a_1a_2-a_3a_6a_1a_4+a_6^2a_1^2+a_5a_3a_6-a_5a_3a_4a_2+a_3^2a_6a_2%
}, \\
z=a_5\frac{a_2a_5^2-a_5a_6a_1-a_5a_4a_3+a_3^2a_6}{%
a_5^2a_2^2-a_5^2a_4+a_5a_4^2a_1-2a_5a_6a_1a_2-a_3a_6a_1a_4+a_6^2a_1^2+a_5a_3a_6-a_5a_3a_4a_2+a_3^2a_6a_2%
}, \\
p=a_6\frac{a_2a_5^2-a_5a_6a_1-a_5a_4a_3+a_3^2a_6}{%
a_5^2a_2^2-a_5^2a_4+a_5a_4^2a_1-2a_5a_6a_1a_2-a_3a_6a_1a_4+a_6^2a_1^2+a_5a_3a_6-a_5a_3a_4a_2+a_3^2a_6a_2%
}.
\end{array}
\right. $
\end{center}

{\bf Proof}\ \ Since $ a(s)$ is Hurwitz stable, Lemma 1 is proved by a direct calculation.

{\bf Lemma {\bf 2}}\ \ Suppose $
s^6+a_1s^5+a_2s^4+a_3s^3+a_4s^2+a_5s+a_6\in H^6,$ then the
following quadratic curve is an  ellipse in  the first quadrant of
the $x$-$y$-$z$-$p$ space:
\begin{center}
$\left\{
\begin{array}{l}
(a_5+a_3y+a_1p-a_2z-a_4x)^2-4(a_2x+z-a_1y-a_3)(a_6x+a_4z-a_3p-a_5y)=0 \\
a_1-x=0 \\
a_5p-a_6z=0
\end{array}
\right. $
\end{center}
and this ellipse is tangent with the line
\begin{center}
$\left\{
\begin{array}{l}
a_1-x=0 \\
a_5p-a_6z=0 \\
a_2x+z-a_1y-a_3=0
\end{array}
\right. $
\end{center}
at
\begin{center}
$\left\{
\begin{array}{l}
x=a_1, \\
y=\frac{-a_5^2+a_5a_4a_1-a_5a_2^2a_1+a_5a_2a_3+a_6a_1^2a_2-a_1a_6a_3}{%
-a_2a_1a_5+a_1^2a_6+a_3a_5}, \\
z=a_5\frac{a_3^2-a_3a_2a_1+a_4a_1^2-a_1a_5}{-a_2a_1a_5+a_1^2a_6+a_3a_5}, \\
p=a_6\frac{a_3^2-a_3a_2a_1+a_4a_1^2-a_1a_5}{-a_2a_1a_5+a_1^2a_6+a_3a_5},
\end{array}
\right. $
\end{center}
tangent with the line
\begin{center}
$\left\{
\begin{array}{l}
a_1-x=0 \\
a_5p-a_6z=0 \\
a_6x+a_4z-a_3p-a_5y=0
\end{array}
\right. $
\end{center}
at
\begin{center}
$\left\{
\begin{array}{l}
x=a_1, \\
y=\frac{a_4^2a_5a_1-a_3a_4a_1a_6-a_5^2a_4+a_3a_5a_6-a_2a_1a_5a_6+a_1^2a_6^2}{%
-a_5^2a_2+a_5a_4a_3-a_3^2a_6+a_5a_6a_1}, \\
z=a_5\frac{a_5a_4a_1-a_1a_6a_3-a_5^2}{-a_5^2a_2+a_5a_4a_3-a_3^2a_6+a_5a_6a_1}%
, \\
p=a_6\frac{a_5a_4a_1-a_1a_6a_3-a_5^2}{-a_5^2a_2+a_5a_4a_3-a_3^2a_6+a_5a_6a_1}%
..
\end{array}
\right. $
\end{center}

{\bf Proof}\ \ Since $ a(s)$ is Hurwitz stable, Lemma 2 is proved by a direct calculation.

{\bf Lemma {\bf 3}}\ \ Suppose $
s^6+a_1s^5+a_2s^4+a_3s^3+a_4s^2+a_5s+a_6\in H^6,$ then the
following quadratic curve is an  ellipse in  the first quadrant of
the $x$-$y$-$z$-$p$ space:
\begin{center}
$\left\{
\begin{array}{l}
(a_6x+a_4z-a_3p-a_5y)^2-4(a_5+a_3y+a_1p-a_2z-a_4x)(a_5p-a_6z)=0 \\
a_1-x=0 \\
a_2x+z-a_1y-a_3=0
\end{array}
\right. $
\end{center}
and this ellipse is tangent with the line
\begin{center}
$\left\{
\begin{array}{l}
a_1-x=0 \\
a_2x+z-a_1y-a_3=0 \\
a_5+a_3y+a_1p-a_2z-a_4x=0
\end{array}
\right. $
\end{center}
at
\begin{center}
$\left\{
\begin{array}{l}
x=a_1, \\
y=-\frac{a_3a_5+a_3a_2^2a_1+a_1^2a_6-a_2a_3^2-a_4a_2a_1^2}{%
a_3^2-a_3a_2a_1+a_4a_1^2-a_1a_5}, \\
z=-\frac{a_2a_1a_3^2-a_2a_1^2a_5+2a_1a_3a_5+a_1^3a_6-a_3^3-a_3a_4a_1^2}{%
a_3^2-a_3a_2a_1+a_4a_1^2-a_1a_5}, \\
p=-\frac{%
2a_5a_4a_1+a_4a_1a_2a_3-a_5^2-a_5a_2^2a_1+a_5a_2a_3-a_1a_6a_3+a_6a_1^2a_2-a_4a_3^2-a_4^2a_1^2%
}{a_3^2-a_3a_2a_1+a_4a_1^2-a_1a_5},
\end{array}
\right. $
\end{center}
tangent with the line
\begin{center}
$\left\{
\begin{array}{l}
a_1-x=0 \\
a_2x+z-a_1y-a_3=0 \\
a_5p-a_6z=0
\end{array}
\right. $
\end{center}
at
\begin{center}
$\left\{
\begin{array}{l}
x=a_1, \\
y=-\frac{a_5a_6a_1-a_5a_4a_2a_1+a_5a_4a_3+a_3a_6a_2a_1-a_3^2a_6}{%
a_5a_4a_1-a_1a_6a_3-a_5^2}, \\
z=-a_5\frac{-a_2a_1a_5+a_1^2a_6+a_3a_5}{a_5a_4a_1-a_1a_6a_3-a_5^2}, \\
p=-a_6\frac{-a_2a_1a_5+a_1^2a_6+a_3a_5}{a_5a_4a_1-a_1a_6a_3-a_5^2}.
\end{array}
\right. $
\end{center}

{\bf Proof}\ \ Since $ a(s)$ is Hurwitz stable, Lemma 3 is proved
by a direct calculation.

For notational simplicity, denote
$$\begin{array}{ll}
\Omega _{e1}^a:=\{(x,y,z,p)|&
(a_2x+z-a_1y-a_3)^2-4(a_1-x)(a_5+a_3y+a_1p-a_2z \\
&-a_4x)<0, a_6x+a_4z-a_3p-a_5y=0, a_5p-a_6z=0 \}
\end{array}
$$ %
$$\begin{array}{ll}
\Omega _{e2}^a:=\{(x,y,z,p)|&
(a_5+a_3y+a_1p-a_2z-a_4x)^2-4(a_2x+z-a_1y-a_3)(a_6x \\
&+a_4z-a_3p-a_5y)<0, a_1-x=0, a_5p-a_6z=0 \}
\end{array}
$$ %
$$\begin{array}{ll}
\Omega _{e3}^a:=\{(x,y,z,p)|&
(a_6x+a_4z-a_3p-a_5y)^2-4(a_5+a_3y+a_1p-a_2z \\
&-a_4x)(a_5p-a_6z)<0, a_1-x=0, a_2x+z-a_1y-a_3=0 \}
\end{array}
$$ %

$$\begin{array}{ll}
\Omega _{e1}^b:=\{(x,y,z,p)|&
(b_2x+z-b_1y-b_3)^2-4(b_1-x)(b_5+b_3y+b_1p-b_2z \\
&-b_4x)<0, b_6x+b_4z-b_3p-b_5y=0, b_5p-b_6z=0 \}
\end{array}
$$ %
$$\begin{array}{ll}
\Omega _{e2}^b:=\{(x,y,z,p)|&
(b_5+b_3y+b_1p-b_2z-b_4x)^2-4(b_2x+z-b_1y-b_3)(b_6x \\
&+b_4z-b_3p-b_5y)<0, b_1-x=0, b_5p-b_6z=0 \}
\end{array}
$$ %
$$\begin{array}{ll}
\Omega _{e3}^b:=\{(x,y,z,p)|&
(b_6x+b_4z-b_3p-b_5y)^2-4(b_5+b_3y+b_1p-b_2z \\
&-b_4x)(b_5p-b_6z)<0, b_1-x=0, b_2x+z-b_1y-b_3=0 \}
\end{array}
$$ %

In what follows, $(A, B)$ stands for the set of points in the line
segment connecting the point $A$ and the point $B$ in the
$x$-$y$-$z$-$p$ space, not including the endpoints $A$ and $B$.
Denote
$$\Omega ^a:=\{(x,y,z,p)|(x,y,z,p)\in (A,B)\cup (A,C)\cup (B,C),
\forall A \in \Omega _{e1}^a,\forall B \in \Omega _{e2}^a,\forall
C \in \Omega _{e3}^a\}$$
$$\Omega ^b:=\{(x,y,z,p)|(x,y,z,p)\in (A,B)\cup (A,C)\cup (B,C),
\forall A \in \Omega _{e1}^b,\forall B \in \Omega _{e2}^b, \forall
C \in \Omega _{e3}^b\}$$

{\bf  Lemma {\bf 4}}\ \ Suppose $
a(s)=s^6+a_1s^5+a_2s^4+a_3s^3+a_4s^2+a_5s+a_6\in H^6,
b(s)=s^6+b_1s^5+b_2s^4+b_3s^3+b_4s^2+b_5s+b_6\in H^6,$ if $ \Omega
^a \cap \Omega ^b\neq \phi,$ take $(x,y,z,p)\in \Omega ^a \cap
\Omega ^b, $ and let $
c(s):=s^5+(x-\varepsilon)s^4+ys^3+zs^2+ps+\varepsilon$ ($
\varepsilon $ is a sufficiently small positive number), then for
$\displaystyle\frac{c(s)}{a(s)}$ and
$\displaystyle\frac{c(s)}{b(s)}$, we have $\forall \omega \in R,
\mbox{Re} [\displaystyle\frac {c(j\omega )}{a(j\omega )}]>0$  and
$ \mbox{Re} [\displaystyle\frac {c(j\omega )}{b(j\omega )}]>0.$

{\bf\bf Proof}\ \ Suppose $(x,y,z,p)\in \Omega ^a \cap \Omega ^b, $ let $%
c(s):=s^5+(x-\varepsilon)s^4+ys^3+zs^2+ps+\varepsilon, \varepsilon
>0,\varepsilon $ sufficiently small.

$\forall \omega \in R,$ consider
$$
\begin{array}{ll}
\mbox{Re}[\displaystyle\frac{c(j\omega )}{a(j\omega )}]= & \displaystyle%
\frac 1{\mid a(j\omega )\mid ^2}[(a_1-x)\omega
^{10}+(a_2x+z-a_1y-a_3)\omega
^8+(a_5+a_3y+a_1p \\
& -a_2z-a_4x)\omega ^6+(a_6x+a_4z-a_3p-a_5y)\omega
^4+(a_5p-a_6z)\omega ^2 \\
& +\varepsilon (\omega ^{10}-a_2\omega ^8+(a_4-1)\omega
^6+(-a_6+a_2)\omega ^4-a_4\omega ^2+a_6)]
\end{array}
$$

In order to prove that $\forall \omega \in R,\mbox{Re}[\displaystyle\frac{%
c(j\omega )}{a(j\omega )}]>0,$ let $t=\omega ^2,$ we only need to
prove that, for any $\varepsilon >0,\varepsilon $ sufficiently
small, the following polynomial $f_1(t)$ satisfies
$$
\begin{array}{ll}
f_1(t):= & t[(a_1-x)t^4+(a_2x+z-a_1y-a_3)t^3+(a_5+a_3y+a_1p \\
& -a_2z-a_4x)t^2+(a_6x+a_4z-a_3p-a_5y)t+(a_5p-a_6z)] \\
& +\varepsilon (t^5-a_2t^4+(a_4-1)t^3+(-a_6+a_2)t^2-a_4t+a_6)]>0,
\ \ \forall t\in [0,+\infty ).
\end{array}
$$

Since $(x,y,z,p)\in \Omega ^a,$ by the definition of $\Omega ^a,$
it is easy to know that
$$
\begin{array}{ll}
g_1(t):= & (a_1-x)t^4+(a_2x+z-a_1y-a_3)t^3+(a_5+a_3y+a_1p-a_2z-a_4x)t^2 \\
& +(a_6x+a_4z-a_3p-a_5y)t+(a_5p-a_6z)>0,\ \ \forall t\in
[0,+\infty ).
\end{array}
$$
Moreover, we obviously have $f_1 (0)>0,$ and for any $\varepsilon >0$, when $%
t$ is a sufficiently large or sufficiently small positive number, we have $%
f_1 (t)>0,$ namely, there exist $0<t_1<t_2$ such that, for all
$\varepsilon
>0$, $t\in [0,t_1]\cup [t_2,+\infty)$, we have $f_1 (t)>0.$

Denote
$$
M_1=\dinf_{t\in [t_1,t_2]}tg_1(t),
$$
$$
N_1=\dsup_{t\in
[t_1,t_2]}|t^5-a_2t^4+(a_4-1)t^3+(-a_6+a_2)t^2-a_4t+a_6|.
$$
Then $M_1>0$ and $N_1>0.$ Choosing $0<\varepsilon <\displaystyle\frac{M_1}{%
N_1},$ by a direct calculation, we have
$$
\begin{array}{ll}
f_1(t):= & t[(a_1-x)t^4+(a_2x+z-a_1y-a_3)t^3+(a_5+a_3y+a_1p \\
& -a_2z-a_4x)t^2+(a_6x+a_4z-a_3p-a_5y)t+(a_5p-a_6z)] \\
& +\varepsilon (t^5-a_2t^4+(a_4-1)t^3+(-a_6+a_2)t^2-a_4t+a_6)>0,\
\ \forall t\in [0,+\infty ).
\end{array}
$$
Namely,
$$
\forall \omega \in R,\mbox{Re}[\displaystyle\frac{c(j\omega )}{a(j\omega )}%
]>0.
$$

Similarly, since $(x,y,z,p)\in \Omega ^b,$ there exist $0<t_3<t_4$
such
that, for all $\varepsilon >0$, $t\in [0,t_3]\cup [t_4,+\infty )$, we have $%
f_2(t)>0,$where
$$
\begin{array}{ll}
f_2(t):= & t[(b_1-x)t^4+(b_2x+z-b_1y-b_3)t^3+(b_5+b_3y+b_1p \\
& -b_2z-b_4x)t^2+(b_6x+b_4z-b_3p-b_5y)t+(b_5p-b_6z)] \\
& +\varepsilon (t^5-b_2t^4+(b_4-1)t^3+(-b_6+b_2)t^2-b_4t+b_6)]
\end{array}
$$

Denote
$$
\begin{array}{ll}
g_2(t):= & (b_1-x)t^4+(b_2x+z-b_1y-b_3)t^3+(b_5+b_3y+b_1p-b_2z-b_4x)t^2 \\
& +(b_6x+b_4z-b_3p-b_5y)t+(b_5p-b_6z),
\end{array}
$$
$$
M_2=\dinf_{t\in [t_3,t_4]}tg_2(t),
$$
$$
N_2=\dsup_{t\in
[t_3,t_4]}|t^5-b_2t^4+(b_4-1)t^3+(-b_6+b_2)t^2-b_4t+b_6|.
$$
Then $M_2>0$ and $N_2>0.$ Choosing $0<\varepsilon <\displaystyle\frac{M_2}{%
N_2},$ we have
$$
\forall \omega \in R,\mbox{Re}[\displaystyle\frac{c(j\omega )}{b(j\omega )}%
]>0.
$$
Thus, by choosing $0<\varepsilon <\min \{{\displaystyle\frac{M_1}{N_1},%
\displaystyle\frac{M_2}{N_2}}\},$ Lemma 4 is proved.

{\bf  Lemma {\bf 5}}\ \ Suppose $
a(s)=s^6+a_1s^5+a_2s^4+a_3s^3+a_4s^2+a_5s+a_6\in H^6,
b(s)=s^6+b_1s^5+b_2s^4+b_3s^3+b_4s^2+b_5s+b_6\in H^6,$ if $\lambda
b(s)+(1-\lambda )a(s)\in H^6,\lambda \in [0,1],$ then $ \Omega ^a
\cap \Omega ^b\neq \phi $

{\bf  Proof }\ \  If $ \forall \lambda \in
[0,1],{a_\lambda}(s):=\lambda b(s)+(1-\lambda )a(s)\in H^6,$ by
Lemmas 1-3, for any $ \lambda \in [0,1],$ $\Omega
_{e1}^{a_\lambda}, \Omega _{e2}^{a_\lambda}$ and $\Omega
_{e3}^{a_\lambda}$ are three ellipses in the first quadrant of the
$x$-$y$-$z$-$p$ space, denote
$$\Omega
^{a_\lambda}:=\{(x,y,z,p)|(x,y,z,p)\in (A,B)\cup (A,C)\cup (B,C),
\forall A \in \Omega _{e1}^{a_\lambda}, \forall B \in \Omega
_{e2}^{a_\lambda},\forall C \in \Omega _{e3}^{a_\lambda}\}.
$$%

Apparently, when $\lambda $ changes continuously from $0$ to $1$,
$\Omega^{a_\lambda}$ will change continuously from $\Omega^{a}$ to
$\Omega^{b},$ $\Omega_{e1}^{a_\lambda}$ will change continuously
from $\Omega_{e1}^{a}$ to
$\Omega_{e1}^{b},\Omega_{e2}^{a_\lambda}$ will change continuously
from $\Omega_{e2}^{a}$ to $\Omega_{e2}^{b},$ and
$\Omega_{e3}^{a_\lambda}$ will change continuously from
$\Omega_{e3}^{a}$ to $\Omega_{e3}^{b}.$

Now assume $ \Omega ^a \cap \Omega ^b= \phi,$ by the definitions
of $ \Omega ^a $ and $ \Omega ^b,$ and Lemmas 1-3, $ \exists u
> 0,$ $v > 0,u\neq a_1, u\neq b_1, $ and $\exists k\in \{1,2,3\},$
such that the following plane in the $x$-$y$-$z$-$p$ space
$$ l:\ \ \  \displaystyle\frac  {x}{u}+\displaystyle\frac  {y}{v}
+\displaystyle\frac  {z}{w}+\displaystyle\frac  {p}{r}=1$$ \\
separates $\Omega^{a}$ and $\Omega^{b},$ meanwhile, $l$ is tangent
with $ \Omega _{e1}^a ,\Omega _{e2}^a, \Omega _{e3}^a$ and $
\Omega _{ek}^b $ simultaneously (or tangent with $ \Omega _{e1}^b
,\Omega _{e2}^b, \Omega _{e3}^b$ and $ \Omega _{ek}^a $
simultaneously ).

Without loss of generality, suppose that $l$ is tangent with $
\Omega _{e1}^a ,\Omega _{e2}^a, \Omega _{e3}^a$ and $ \Omega
_{ek}^b $ simultaneously.

Since $ l$ is tangent with $ \Omega _{e1}^a, \Omega _{e2}^a $ and
$ \Omega _{e3}^a $ simultaneously, $ a(s)$ is Hurwitz stable and
$u > 0,u\neq a_1, v > 0, $ by a lengthy calculation, we get that
the necessary and sufficient condition for $ l$ being tangent with
$ \Omega _{e1}^a, \Omega _{e2}^a $ and $ \Omega _{e3}^a $
simultaneously is
\begin{equation}
uv^3-a_1v^3-a_2uv^2+a_3v^2+a_4uv-a_5v-a_6u=0,w=-uv,r=-v^2
\end {equation}

Since $w=-uv,r=-v^2, l$ is tangent with $ \Omega _{ek}^b $, by a
direct calculation, we have
\begin{equation}
uv^3-b_1v^3-b_2uv^2+b_3v^2+b_4uv-b_5v-b_6u=0
\end {equation}\\%

>From $ (1)$  and $ (2),$ we obviously have $\forall \lambda \in [0,1],$
\begin{equation}
uv^3-a_{\lambda 1}v^3-a_{\lambda 2}uv^2+a_{\lambda
3}v^2+a_{\lambda 4}uv-a_{\lambda 5}v-a_{\lambda 6}u=0,w=-uv,r=-v^2
\end {equation}\\
where $a_{\lambda i}:=a_i+\lambda (b_i-a_i),i=1,2,3,4,5,6.$ $(3)$
shows that $ l$ is also tangent with $ \Omega _{ek}^{a_{\lambda
}}( \forall \lambda \in [0,1])$. But $l $ separates $ \Omega
_{ek}^a $ and $ \Omega _{ek}^b, $ and when $\lambda $ changes
continuously from $0$ to $1$, $\Omega_{ek}^{a_\lambda}$ will
change continuously from $\Omega_{ek}^{a}$ to $\Omega_{ek}^{b},$
which is obviously impossible. This completes the proof.

From Theorem 2.4 in \cite{WY00}, or the proof of Lemma 5 in
\cite{YH99}, we have

{\bf  Lemma {\bf 6}}\ \ Suppose $
a(s)=s^6+a_1s^5+a_2s^4+a_3s^3+a_4s^2+a_5s+a_6\in H^6,
b(s)=s^6+b_1s^5+b_2s^4+b_3s^3+b_4s^2+b_5s+b_6\in H^6,
c(s)=s^5+xs^4+ys^3+zs^2+ps+q,$ if $ \forall \omega \in R,
\mbox{Re} [\displaystyle\frac {c(j\omega )}{a(j\omega )}]>0$  and
$ \mbox{Re} [\displaystyle\frac {c(j\omega )}{b(j\omega )}]>0,$
take
\begin{center}
$ \stackrel{\sim }{c}(s):=c(s)+\delta \cdot d(s),\ \ \delta
>0,\ \delta $ sufficiently small
\end{center}
(where $d(s)$ is an arbitrarily given monic sixth-order
polynomial), then $ \displaystyle\frac {\stackrel{\sim
}{c}(s)}{a(s)}$  and
$ \displaystyle\frac {\stackrel{\sim }{c}(s)}{%
b(s)}$ are both strictly positive real.

The sufficiency of Theorem 1 is now proved by combining Lemmas 1-6. %

{\bf  Remark {\bf 1}} \ \ From the proof of Theorem 1, we can see
that this paper not only proves the existence,
but also provides a design method. %

{\bf  Remark {\bf 2}} \ \ The method provided in this paper is
constructive, and is insightful and helpful in solving the general
robust SPR synthesis problem. This subject is currently under
investigation.%

{\bf  Remark {\bf 3}} \ \
Our results can easily be generalized to discrete-time case. %

{\bf  Remark {\bf 4}} \ \
If $\displaystyle\frac{c(s)}{a(s)}$ and $\displaystyle\frac{c(s)}{b(s)}$ are both SPR,
it is easy to know that
$\forall \lambda \in [0, 1], \displaystyle \frac{c(s)}{\lambda a(s) + (1- \lambda) b(s)}$ is also SPR.

{\bf  Remark {\bf 5}} \ \ The stability of polynomial segment can
be checked by many efficient methods, e.g., eigenvalue method,
root locus method, value set method, etc. \cite{Bar94,BCK95}.

}}

\vskip 20pt
\vspace*{1\baselineskip}

\end{document}